\documentclass[11pt]{article}

\usepackage{graphicx,latexsym,amssymb}

\newtheorem{lemma}{{\bf Lemma}}
\newtheorem{theo}{{\bf Theorem}}
\newtheorem{prop}{{\bf Proposition}}
\newtheorem{coro}{{\bf Corollary}}

\newtheorem{remark}{{\bf Remark}}

\font\bbb=msbm9.5 scaled\magstep1

\newcommand{\CC}{\mbox{\bbb C}}
\newcommand{\FF}{\mbox{\bbb F}}
\newcommand{\HH}{\mbox{\bbb H}}

\newcommand{\ZZ}{\mbox{\bbb Z}}

\newcommand{\TPSS}{S^{\hspace{.2mm}3}\rotatebox{90}{\ensuremath{\ltimes}} S^{\hspace{.1mm}1}}
\newcommand{\TPSSD}{S^{\hspace{.2mm}d-1}\rotatebox{90}{\ensuremath{\ltimes}} S^{\hspace{.1mm}1}}

\newcommand{\XB}{X \hspace{-2.7mm}^{^{\mbox{\bf --}}}}

 \newcommand{\TPPSS}{\kern.24em
 \rule width.08em
   height1.5ex depth-.08ex
 \kern-.36em \times}

\begin{document}

\title{\bf On Walkup's class \boldmath{${\cal K}(d)$} and a minimal triangulation
of \boldmath{$(\TPSS)^{\#3}$}}
\author{{\bf Bhaskar Bagchi}, {\bf Basudeb Datta}
}

\date{}

\maketitle

\vspace{-5mm}

\noindent {\small $^{\rm a}$ Theoretical Statistics and
Mathematics Unit, Indian Statistical Institute,  Bangalore
560\,059, India.

\smallskip

\noindent $^{\rm b}$ Department of Mathematics, Indian Institute
of Science, Bangalore 560\,012,  India.}\footnotetext{{\em E-mail
addresses:} bbagchi@isibang.ac.in (B. Bagchi),
dattab@math.iisc.ernet.in (B. Datta).}

\begin{center}
%\date{Revised\,: February 26, 2011}
\date{To appear in `Discrete Mathematics'}
\end{center}

\hrule

\bigskip

\noindent {\bf Abstract}

{\small For $d \geq 2$, Walkup's class ${\cal K}(d)$ consists of
the $d$-dimensional simplicial complexes all whose vertex-links
are stacked $(d-1)$-spheres. Kalai showed that for $d\geq 4$, all
connected members of ${\cal K}(d)$ are obtained from stacked
$d$-spheres by finitely many elementary handle additions.
According to a result of Walkup, the face vector of any
triangulated 4-manifold $X$ with Euler characteristic $\chi$
satisfies $f_1 \geq 5f_0 - \frac{15}{2} \chi$, with equality only
for $X \in {\cal K}(4)$. K\"{u}hnel observed that this implies
$f_0(f_0 - 11) \geq -15\chi$, with equality only for 2-neighborly
members of ${\cal K}(4)$. K\"{u}hnel also asked if there is a
triangulated 4-manifold with $f_0 = 15$, $\chi = -4$ (attaining
equality in his lower bound). In this paper, guided by Kalai's
theorem, we show that indeed there is such a triangulation. It
triangulates the connected sum of three copies of the twisted
sphere product $\TPSS$. Because of K\"{u}hnel's inequality, the
given triangulation of this manifold is a vertex-minimal
triangulation. By a recent result of Effenberger, the
triangulation constructed here is tight. Apart from the neighborly
2-manifolds and the infinite family of $(2d+ 3)$-vertex sphere
products $S^{\, d-1} \times S^1$ (twisted for $d$ odd), only
fourteen tight triangulated manifolds were known so far. The
present construction yields a new member of this sporadic family.
We also present a self-contained proof of Kalai's result. }

\bigskip

{\small {\em MSC 2000\,:} 57Q15, 57R05.

{\em Keywords:} Stacked spheres; Triangulated manifolds; Tight
triangulations.

}

\bigskip

\hrule

\section{Walkup's class \boldmath{${\cal K}(d)$}}

A {\em weak pseudomanifold} without (respectively, with) boundary
is a pure simplicial complex in which each face of co-dimension
one is in exactly (respectively, at most) two facets (face of
maximum dimension). The {\em dual graph} $\Lambda (X)$ of a weak
pseudomanifold $X$ is the graph (simplicial complex of dimension
$\leq 1$) whose vertices are the facets of $X$, two such vertices
being adjacent in $\Lambda (X)$ if the corresponding facets of $X$
meet in a co-dimension one face. We say that $X$ is a {\em
pseudomanifold} if $\Lambda(X)$ is connected. Any triangulation
of a closed and connected manifold is automatically a
pseudomanifold without boundary.

For a simplicial complex $X$ of dimension $d$, $f_j = f_j(X)$
denotes the number of $j$-dimensional faces of $X$ ($0\leq j\leq
d$), and the vector $f(X) = (f_0, \dots, f_d)$ is called the {\em
face vector} of $X$.

A {\em stacked ball} of dimension $d$ (in short, a stacked
$d$-ball) may be defined as a $d$-dimensional pseudomanifold $X$
with boundary such that $\Lambda(X)$ is a tree. (We recall that a
tree is a minimally connected graph, i.e., a connected graph
which is disconnected by the removal of any of its edges.) A {\em
stacked $d$-sphere} may be defined as the boundary of a stacked
$(d+1)$-ball. Since a tree on at least two vertices has (at least
two) end vertices, a trivial induction shows that a stacked
$d$-ball actually triangulates a topological $d$-ball, and hence
a stacked $d$-sphere triangulates a topological $d$-sphere.
By the same reason, a simplicial complex is a stacked $d$-sphere
if and only if it is obtained from the standard sphere $S^{\,d}_{d
+2}$ by finite sequence of starring vertices in facets (see
\cite[Proposition 4]{bd13}).  (This last condition is usually
used to define stacked spheres in the literature.) Since an
$n$-vertex stacked $d$-sphere is obtained from $S^{\,d}_{d + 2}$
by $(n - d - 2)$ starring and each starring induces ${d+1 \choose
j}$ new $j$-faces and retains all the old $j$-faces for $1\leq j
< d$ (respectively, kills only one old $j$-face for $j=d$), it
follows that it has $(n-d-2){d+1 \choose j} + {d+2 \choose j+1}$
$j$-faces for $1\leq j < d$, and $(n-d-2)d + (d+2)$ facets. On
simplifying, we get\,:

\begin{lemma}$\!\!\!${\bf .} \label{L1}
The face vector of any $d$-dimensional stacked sphere satisfies
$$
f_j = \left\{\begin{array}{ll}
 {d+1 \choose j}f_0 - j{d+2 \choose j+1}, & \mbox{if } ~~ 1\leq j<d \\
 df_0 - (d+2)(d-1), & \mbox{if } ~~ j=d.
 \end{array}
 \right.
$$
\end{lemma}

In \cite{wa}, Walkup defined the class ${\cal K}(d)$ as the
family of all $d$-dimensional simplicial complexes all whose
vertex-links are stacked $(d - 1)$-spheres. Clearly, all the
members of ${\cal K}(d)$ are triangulated manifolds and, indeed,
for $d\leq 2$, ${\cal K}(d)$ consists of all triangulated
$d$-manifolds.

\begin{prop}$\!\!\!${\bf .} \label{P1}
Let $d$ be an even number and let $X$ be a connected member of
${\cal K}(d)$ with Euler characteristic $\chi$. Then the
face vector of $X$ is given by
$$
f_j = \left\{\begin{array}{ll}
 {d+1 \choose j}f_0 - \frac{j}{2}{d+2 \choose j+1}\chi,
 & \mbox{ if } ~~ 1\leq j<d \\
 df_0 - \frac{1}{2}(d+2)(d-1)\chi, & \mbox{if } ~~ j=d.
 \end{array}
 \right.
$$
\end{prop}

\noindent {\bf Proof.} Let's count in two ways the number of
ordered pairs $(x, \tau)$, where $\tau$ is a $j$-face of $X$ and
$x\in \tau$ is a vertex. This yields the formula
$$
f_j = \frac{1}{j+1}\sum_{x\in V(X)} f_{j-1}({\rm lk}(x)).
$$
Let, as usual, $\deg(x)$ denote the degree of $x$ in $X$ (i.e.,
the number of vertices in ${\rm lk}(x)$). Since all the
vertex-links ${\rm lk}(x)$ of $X$ are stacked $(d-1)$-spheres,
Lemma \ref{L1} applied to these links shows that

$$
f_j = \left\{\begin{array}{ll}
 \frac{1}{j+1}\displaystyle{\sum_{x\in V(X)}} \left({d \choose j-1}
 \deg(x) - (j-1){d  +1 \choose j}\right), &  ~~ 1\leq j<d \\
 \frac{1}{j+1}\displaystyle{\sum_{x\in V(X)}}((d -1)\deg(x) -
 (d-2)(d+1)), & ~~ j=d.
 \end{array}
 \right.
$$
But $\sum_{x\in V(X)}\deg(x) = 2f_1$. Therefore, we obtain
\begin{equation} \label{E1}
f_j = \left\{\begin{array}{ll}
 \frac{2}{j+1}{d \choose j-1}f_1 - \frac{j-1}{j+1}{d  +1 \choose j}f_0, &  ~~ 1\leq j<d
 \\[2mm]
 \frac{2d-2}{d+1}f_1 - (d-2)f_0, & ~~ j=d.
 \end{array}
 \right.
\end{equation}

Substituting (\ref{E1}) into $\chi = \sum_{j=0}^d (-1)^jf_j$, and
remembering that $d$ is even, we get $\chi = 2af_1 - bf_0$, where
$a := \frac{d-1}{d+1} + \sum_{j=1}^{d-1}(-1)^j\frac{1}{j+1}{d
\choose j-1}$ and $b := d-2 +
\sum_{j=0}^{d-1}(-1)^j\frac{j-1}{j+1}{d+1 \choose j}$. But the
binomial theorem together with Euler's formula, relating his Beta
and Gamma integrals, yields:
$$
\displaystyle{\sum_{j=1}^{d+1}} (-1)^j\frac{1}{j+1}{d\choose j-1}
= -\int_0^1 (1-x)^d x dx = -\beta(2, d+1) = -\frac{1}{(d+1)(d+
2)},
$$
and
$$
\displaystyle{\sum_{j=0}^{d+1}} (-1)^j\frac{1}{j+1}{d+1\choose j}
= \int_0^1 (1-x)^{d+1} dx = \frac{1}{d+ 2}.
$$
Hence (still remembering that $d$ is even), we get $a = 2/(d+2)$
and $b = 4/(d+1)(d+2)$. Thus, $\chi = 4f_0/(d+2) - 4f_1/(d+ 1)(d+
2)$. In other words, $f_1 = (d+1)f_0 - \frac{1}{2}{d +2\choose 2}
\chi$. Substituting this value of $f_1$ in (\ref{E1}), we get the
expression for $f_j$ in terms of $f_0$ and $\chi$, as claimed.
\hfill $\Box$

\medskip

Notice that, till the proof of (\ref{E1}), we have not used the
assumption that $d$ is even. Thus (\ref{E1}) is valid for all
dimensions $d$. However, there is no further simplification when
$d$ is odd.

A simplicial complex is said to be {\em $2$-neighborly} if any two
vertices are joined by an edge, i.e., $f_1 = {f_0 \choose 2}$.
Thus Proposition \ref{P1} has the following immediate consequence:

\begin{coro}$\!\!\!${\bf .} \label{C1}
Let $d$ be an even number and $X$ be a connected member of ${\cal
K}(d)$ with Euler characteristic $\chi$. Then the face vector of
$X$ satisfies $f_0(f_0 - 2d - 3) \geq - {d+2 \choose 2}\chi$, and
equality holds if and only if $X$ is $2$-neighborly.
\end{coro}

\noindent {\bf Proof.} This is immediate on substituting $f_1 =
(d+1)f_0 - \frac12{d+2 \choose 2}\chi$ in the trivial inequality
$f_1 \leq {f_0\choose 2}$.  \hfill $\Box$

\begin{remark}$\!\!\!${\bf .} \label{R1}
{\rm Let $X$ be a connected member of ${\cal K}(d)$, $d\geq 4$
even, and $\FF$ be a field such that $X$ is $\FF$-orientable. Let
$\beta_i = \beta_i(X ; \FF)$ be the corresponding Betti numbers.
Then Kalai's theorem (Proposition \ref{P3} below) implies that
the Euler characteristic $\chi$ of $X$ is given by $\chi = 2 -
2\beta_1$. Therefore, the inequality of Corollary \ref{C1} may
be rewritten as\,:
$$
{f_0 - d - 1 \choose 2} \geq {d+2 \choose 2}\beta_1.
$$
In \cite{lss}, Lutz, Sulanke and Swartz have shown that this last
inequality holds for all $\FF$-orientable connected triangulated
$d$-manifolds $M$, for $d\geq 3$. Further, equality holds here if and
only if $M$ is a 2-neighbourly member of ${\cal K}(d)$. They call
a $d$-manifold {\em tight neighbourly} if it attains equality in
their bound. Thus, tight neighbourly $d$-manifolds are precisely
the 2-neighbourly members of ${\cal K}(d)$. For instance, the
triangulated 4-manifold $M^{4}_{15}$ of Section 2 below is tight
neighbourly.}
\end{remark}

\begin{prop} {\bf (Walkup \cite{wa}, K\"{u}hnel
\cite{ku}).} \label{P2} Let $X$ be a connected triangulated
$4$-manifold with Euler characteristic $\chi$. Then the
face vector of $X$ satisfies the following.
\begin{enumerate} \vspace{-2mm}
\item[{\rm (a)}]
$$
f_j \geq \left\{\begin{array}{ll}
 {5 \choose j}f_0 - \frac{j}{2}{6 \choose j+1}\chi,
 & \mbox{ if } ~~ 1\leq j<4 \\
 4f_0 - 9\chi, & \mbox{if } ~~ j=4.
 \end{array}
 \right.
$$
Further, equality holds here for some $j \geq 1$ if and only if $X
\in {\cal K}(4)$. \item[{\rm (b)}] $f_0(f_0-11) \geq - 15  \chi$,
and equality holds here if and only if $X$ is a $2$-neighborly
member of ${\cal K}(4)$.
\end{enumerate}
\end{prop}

\noindent {\bf Proof.} As a well-known consequence of the
Dehn-Sommerville equations, the face vector of $X$ satisfies (cf.
\cite{ku}) $f_2 = 4f_1 - 10(f_0 - \chi)$, $f_3 = 5f_1 - 15(f_0 -
\chi)$ and $f_4 = 2f_1 - 6(f_0 - \chi)$. Therefore, to prove Part
(a), it suffices to do the case $j =1$: $f_1 \geq 5 f_0 -
15\chi/2$, with equality only for $X \in {\cal K}(4)$. But,
applying the lower bound theorem (LBT) for normal pseudomanifolds
(cf. \cite{bd9}) to the vertex links of $X$, we get $f_2 =
\frac13\sum_{x\in V(X)} f_1({\rm lk}(x)) \geq \frac13 \sum_{x\in
V(X)} (4\deg(x) - 10) = \frac83 f_1 -\frac{10}{3}f_0$. On
substituting $f_2 = 4f_1 - 10(f_0 - \chi)$, this simplifies to
$f_1 \geq 5f_0 - \frac{15}{2}\chi$. Since equality in the LBT
holds only for stacked spheres, equality holds only for $X\in
{\cal K}(4)$. This proves (a).

In conjunction with the trivial inequality $f_1 \leq {f_0 \choose
2}$, Part (a) implies Part (b). \hfill $\Box$

\medskip

Clearly, any $d$-dimensional weak pseudomanifold without boundary
has at least $d+2$ vertices, with equality only for the standard
$d$-sphere $S^{\,d}_{d+2}$ (the boundary complex of a $(d+
1)$-simplex). $S^{\,d}_{d+2}$ is a stacked
$d$-sphere: it is the boundary of the standard $(d+1)$-ball
$B^{\,d+1}_{d+2}$ (the face complex of a $(d+1)$-simplex). Since any tree on at
least two vertices has at least two end vertices (i.e., vertices
of degree one), the following lemma is immediate from the
definitions of stacked balls and stacked spheres. (See \cite{bd9}
for a proof.)

\begin{lemma}$\!\!\!${\bf .} \label{L2}
Let $X$ be a stacked $d$-sphere. \begin{enumerate} \vspace{-2mm}
\item[{\rm (a)}] Then $X$ has at least two vertices of
$($minimum$)$ degree $d+1$. \vspace{-2mm} \item[{\rm (b)}] Let $X$
have $f_0 > d+2$ vertices. Suppose $x$ is a vertex of degree
$d+1$. Let $\sigma$ denote the set of neighbors of $x$ in $X$. Let
$X_0$ be the pure simplicial complex whose facets are $\sigma$
together with the facets of $X$ not containing $x$. Then $X_0$ is
a stacked $d$-sphere.
\end{enumerate}
\end{lemma}

\begin{lemma}$\!\!\!${\bf .} \label{L3}
Let $X$ be a stacked sphere of dimension $d\geq 2$ with edge graph
$(1$-skeleton$)$ $G$. Let $\XB$ be the simplicial complex whose
faces are all the cliques $($sets of mutually adjacent vertices$)$
of $G$. Then $\XB$ is a stacked $(d+ 1)$-ball whose boundary is
$X$.
\end{lemma}

\noindent {\bf Proof.} Let $X$ have $n$ vertices. If $n = d + 2$
then $X = S^{\,d}_{d + 2}$ and $\XB = B^{\,d + 1}_{d + 2}$, and
there is nothing to prove. So assume that $n > d + 2$ and we have
the result for all stacked $d$-spheres with fewer vertices. Let
$x$, $\sigma$, $X_0$ be as in Lemma 2 (b). Notice that (as $d \geq
2$), the edge graph $G_0$ of the $(n-1)$-vertex stacked $d$-sphere
$X_0$ is obtained from $G$ by deleting all the edges through $x$
(and the vertex $x$ itself). Therefore, the cliques of $G$ are
$\alpha \cup \{x\}$, where $\alpha \subseteq \sigma$; and the
cliques of $G_0$. Hence the facets of $\XB$ are $\tilde{\sigma}
:= \sigma \cup \{x\}$ and the facets of the stacked $(d + 1)$-ball
${\XB}_0$. Thus the dual graph $\Lambda(\XB)$ is obtained from
the tree $\Lambda({\XB}_0)$ by adding an end vertex
($\tilde{\sigma}$). So, $\Lambda(\XB)$ is a tree, i.e., $\XB$ is
a stacked $(d + 1)$-ball. Since $X_0$ is the boundary of
${\XB}_0$, it is immediate that $X$ is the boundary of $\XB$.
\hfill $\Box$

\medskip

Notice that Lemma \ref{L3} shows that any stacked sphere is
uniquely determined by its 1-skeleton. (This is, of course,
trivial for $d=1$.)

Now, let $X$ be a member of ${\cal K}(d)$, $d\geq 3$. Let $S$ be
a set of $d+1$ vertices of $X$ such that the induced subcomplex
$X[S]$ of $X$ on the vertex set $S$ is isomorphic to the standard
$(d-1)$-sphere $S^{\,d-1}_{d+1}$. Then, with notations as in
Lemma \ref{L3}, it is clear that for any $x\in S$, $\sigma := S
\setminus \{x\}$ is a $(d-1)$-face in the interior of the stacked
$d$-ball $\overline{{\rm lk}(x})$. The proof of the following
lemma shows that when $d\geq 4$, the converse is also true: if
$\sigma$ is an interior $(d-1)$-face of $\overline{{\rm lk}(x})$
for some vertex $x$, then $X \in {\cal K}(d)$ induces an
$S^{\,d-1}_{d+1}$ on the vertex set $\sigma\cup\{x\}$.

\begin{lemma}$\!\!\!${\bf .} \label{L4} For $d \geq 4$, every member of
${\cal K}(d)$, excepting $S^{\,d}_{d + 2}$, has an $S^{\,d - 1}_{d
+ 1}$ as an induced subcomplex.
\end{lemma}

\noindent {\bf Proof.} Let $X \in {\cal K}(d)$, $X\neq S^{\,d}_{d
+ 2}$. Then $X$ has a vertex of degree $\geq d + 2$. Fix such a
vertex $x$. Then the stacked $d$-ball $\overline{{\rm lk}(x})$
given by Lemma \ref{L3} has an interior $(d - 1)$-face $\sigma$.
(If there was no such $(d - 1)$-face, then we would have
$\overline{{\rm lk}(x}) = B^{\,d}_{d + 1}$, and hence $\deg(x) =
d+1$, contrary to the choice of $x$.) We claim that $X$ induces
an $S^{\,d - 1}_{d + 1}$ on $S := \sigma \cup \{x\}$. Clearly,
every proper subset of $S$, with the possible exception of
$\sigma$, is a face of $X$, while $S$ itself is not a face of $X$
since $\sigma$ is not a boundary face of $\overline{{\rm lk}(x})$.
Therefore, to prove the claim, we need to show that $\sigma\in X$.
Notice that ${\rm lk}(x)$ and $\overline{{\rm lk}(x})$ have the
same $(d-2)$-skeleton. In particular, as $d - 2 \geq 2$ and
$\sigma \in \overline{{\rm lk}(x})$, it follows that each
3-subset of $\sigma$ is in ${\rm lk}(x)$. Therefore, for any
vertex $y \in \sigma$, each 3-subset of $\sigma \cup \{x\}
\setminus \{y\}$ containing $x$ is in ${\rm lk}(y)$. Hence each
2-subset of $\sigma \setminus \{y\}$ is in ${\rm lk}(y)$, i.e.,
$\sigma \setminus \{y\}$ is a clique in the edge graph of ${\rm
lk}(y)$. Hence $\sigma \setminus \{y\} \in \overline{{\rm
lk}(y})$. Since $\sigma \setminus \{y\}$ is a $(d-2)$-face of
$\overline{{\rm lk}(y})$, and ${\rm lk}(y)$ has the same
$(d-2)$-skeleton as $\overline{{\rm lk}(y})$, it follows that
$\sigma \setminus \{y\} \in {\rm lk}(y)$, i.e., $\sigma \in X$, as
was to be shown. \hfill $\Box$

\medskip

Now, let $X$ be a triangulated closed $d$-manifold and
$\sigma_1$, $\sigma_2$ be two facets of $X$. A bijection $\psi :
\sigma_1 \to \sigma_2$ is said to be {\em admissible} if, for
each vertex $x \in \sigma_1$, $x$ and $\psi(x)$ are at distance at
least three in the edge graph of $X$ (i.e., there is no path of
length at most two joining $x$ to $\psi(x)$). In this case, the
triangulated $d$-manifold $X^{\psi}$, obtained from $X \setminus
\{\sigma_1, \sigma_2\}$ by identifying $x$ with $\psi(x)$ for
each $x \in \sigma_1$, is said to be obtained from $X$ by an {\em
elementary handle addition}. Notice that the induced subcomplex of
$X^{\psi}$ on the vertex set $\sigma_1$ ($\approx \sigma_2$) is an
$S^{\,d-1}_{d+1}$. In case $X = X_1 \sqcup X_2$, for
vertex-disjoint subcomplexes $X_1$, $X_2$ of $X$, and $\sigma_1
\in X_1$, $\sigma_2 \in X_2$, any bijection $\psi \colon \sigma_1
\to \sigma_2$ is admissible. In this situation, we write $X_1 \#
X_2$ for $X^{\psi}$, and $X_1 \# X_2$ is called a {\em
$($combinatorial\,$)$ connected sum} of $X_1$  and $X_2$.

In Lemma 1.3 of \cite{bd8}, we have shown (in particular) that if
$Y$ is a connected triangulated closed manifold of dimension $d
\geq 3$, with a vertex set $S$ on which $Y$ induces an $S^{\,d -
1}_{d + 1}$, the above construction can be reversed. Namely, then
there exists a unique triangulated closed $d$-manifold
$\widetilde{Y}$, together with an admissible map $\psi \colon
\sigma_1 \to \sigma_2$, such that $Y = (\widetilde{Y})^{\psi}$,
and $S = \sigma_1 \approx \sigma_2$. The manifold $\widetilde{Y}$
is said to be obtained from $Y$ by a ({\em combinatorial}\,) {\em
handle deletion}. Either $\widetilde{Y}$ is connected, in which
case the first Betti numbers satisfy $\beta_1(Y) =
\beta_1(\widetilde{Y}) + 1$, or else $\widetilde{Y}$ has exactly
two connected components, say $Y_1$ and $Y_2$, and we have $Y =
Y_1 \# Y_2$, in the latter case. It is also easy to see that $Y
\in {\cal K}(d)$ if and only if $\widetilde{Y} \in {\cal K}(d)$
(cf. Lemma 2.6 in \cite{bd8}). We use these results in the
following proof.

\begin{prop} {\bf (Kalai \cite{ka}).} \label{P3} For $d\geq 4$,
a connected simplicial complex $X$ is in ${\cal K}(d)$ if and only
if $X$ is obtained from a stacked $d$-sphere by $\beta_1(X)$
combinatorial handle additions. In consequence, any such $X$
triangulates either $(S^{\,d -1}\!\times S^1)^{\#\beta_1}$ or
$(\TPSSD)^{\#\beta_1}$ according as $X$ is orientable or not.
$($Here $\beta_1 = \beta_1(X).)$
\end{prop}

\noindent {\bf Proof.}  Clearly, stacked $d$-spheres are in
${\cal K}(d)$. Hence so are simplicial complexes obtained from
stacked $d$-spheres by finitely many elementary handle additions.
This proves the ``if" part. We prove the ``only if" part by
induction on the integral Betti number $\beta_1(X)$. To start the
induction, we need\,:

\smallskip

\noindent {\em Claim\,:}  For $d\geq 4$, $X\in {\cal K}(d)$ and
$\beta_1(X) = 0$ imply $X$ is a stacked sphere.

\smallskip

We prove the claim by induction on the number $n$ of vertices in
$X$. If $n  = d+2$ then $X = S^d_{d+2}$, and the result is
obvious. So, assume $n > d+2$ and we have the result for members
of ${\cal K}(d)$ with fewer vertices and vanishing first Betti
number.  By Lemma \ref{L4}, $X$ has an induced subcomplex
isomorphic to $S^{\,d-1}_{d+1}$. Therefore, we may obtain
$\widetilde{X} \in {\cal K}(d)$ by a handle deletion. Then
$\widetilde{X}$ must be disconnected since otherwise we get the
contradiction $\beta_1(X) > \beta(\widetilde{X}) \geq 0$.
Therefore $X = X_1 \# X_2$, where $X_1, X_2 \in {\cal K}(d)$ are
the connected components of $\widetilde{X}$. Since $\beta_1(X_1)
= 0 = \beta_1(X_2)$, induction hypothesis yields that $X_1$,
$X_2$ are both stacked spheres. But the combinatorial connected
sum of stacked spheres is easily seen to be a stacked sphere (cf.
Lemma 2.5 in \cite{bd8}). So, $X$ is a stacked sphere. This
completes the induction, proving the claim.

Thus, we have the ``only if\," part when the Betti number is $0$.
So, assume that the Betti number $\beta_1>0$ and we have the
result for members of ${\cal K}(d)$ with smaller first Betti
number.

If possible, assume that the result is not true, i.e., there
exists a member of ${\cal K}(d)$ with Betti number $\beta_1 >0 $
which can't be obtained from a stacked $d$-sphere by $\beta_1$
combinatorial handle additions. Choose one such member, say $X$,
of ${\cal K}(d)$ with the smallest number of vertices. As before,
obtain $\widetilde{X}$ from $X$ by an combinatorial handle
deletion. If $\widetilde{X}$ is connected then
$\beta_1(\widetilde{X}) = \beta_1 -1$. So, by induction
hypothesis, $\widetilde{X}$ is obtained from a stacked sphere by
$\beta_1(\widetilde{X})$ combinatorial handle additions. Then $X$
is obtained from the same stacked sphere by $\beta_1 =
\beta_1(\widetilde{X}) + 1$ combinatorial handle additions.
Therefore, from our hypothesis, $\widetilde{X}$ is not connected.
So, $\widetilde{X} = X_1 \sqcup X_2$ and $X = X_1\# X_2$, for
some $X_1, X_2 \in {\cal K}(d)$. Then $\beta_1 = \beta_1(X_1) +
\beta_1(X_2)$ and $\beta_1(X_1), \beta_1(X_2) \geq 0$. If
$\beta_1(X_1), \beta_1(X_2) < \beta_1$, then, by induction
hypothesis, $X_i$ is obtained from a stacked sphere $S_i$ by
$\beta_1(X_i)$ combinatorial handle additions, for $1\leq i\leq
2$ and hence $X$ is obtained from the stacked sphere $S_1\# S_2$
by $\beta_1 = \beta_1(X_1) + \beta_1(X_2)$ combinatorial handle
additions. By our assumption, this is not possible. So, one of
$\beta_1(X_1), \beta_1(X_2)$ is equal to $\beta_1$ and the other
is $0$. Assume, without loss, that $\beta_1(X_1) = \beta_1$. This
is a  contradiction to our choice of $X$, since $f_0(X_1) \leq
f_0(X) -1$. Thus, the result  is true for Betti number
$\beta_1$. The first statement now follows by induction.

If $\beta_1 = 1$ then $|X|$ is an $S^{\hspace{.2mm}d - 1}$-bundle
over $S^{\hspace{.2mm}1}$ and hence homeomorphic to $S^{\,d -
1}\!\times S^1$ (if orientable) or $\TPSSD$ (if non-orientable)
(cf. \cite[pages 134--135]{st}). Since $(S^{\,d-1}\!\times S^1) \#
(\TPSSD)$ is homeomorphic to $(\TPSSD) \# (\TPSSD)$, it follows that
$|X|$ is homeomorphic to $(S^{\,d-1}\!\times S^1)^{\#\beta_1}$ (if
orientable) or $(\TPSSD)^{\#\beta_1}$ (if non-orientable). \hfill
$\Box$

\bigskip

We recall that a combinatorial manifold $X$ is said to be {\em
tight} if for every induced subcomplex $Y$ of $X$, the morphism
$H_{\ast}(Y; \ZZ_2) \to H_{\ast}(X; \ZZ_2)$ (induced by the
inclusion map $Y \hookrightarrow X$) is injective (cf. \cite{ku}).
Recently, Effenberger proved\,:

\begin{prop} {\bf (Effenberger \cite{ef}).} \label{P4}
Every $2$-neighborly member of ${\cal K}(d)$ is tight for $d \geq 4$.
\end{prop}

\section{Results}

By Proposition \ref{P2}, any $n$-vertex triangulated connected
4-manifold $X$, with Euler characteristic $\chi$, satisfies $n(n -
11) \geq - 15\chi$. Thus, when $n(n-11) = -15\chi$, $X$ must be a
minimal triangulation of its geometric carrier (requiring the
fewest possible vertices). The smallest values of $n$ for which
equality may hold is $n = 11$. Indeed, there is a unique 11-vertex
4-manifold with $\chi = 0$ (cf \cite{bd8}): it triangulates
$S^{\,3}\!\times S^1$. In \cite{ku}, K\"{u}hnel asked if the next
feasible case $n = 15$, $\chi = -4$ can be realized. Notice that
by Proposition \ref{P2}, any 15-vertex triangulated 4-manifold
with $\chi = -4$ must be a (2-neighborly) member of ${\cal K}(4)$.
By Proposition \ref{P3}, it must arise from a 30-vertex stacked
4-sphere by three elementary handle additions (since it must have
$\beta_1 = 3$). Now, three such operations require three pairs of
facets (each containing five vertices) in the original stacked
sphere, with admissible bijection within each pair. As $30 = 5
\times 6$, it seems reasonable to demand that these six facets in
the sought after 30-vertex stacked 4-sphere be pairwise disjoint,
covering the vertex set (though we are unable to prove that this
must be the case). This strategy works!

\medskip

\noindent {\bf The Construction\,:} Let $B^5_{30}$ denote the
pure 5-dimensional simplicial complex with thirty vertices $a_i$,
$a_i^{\prime}$, $b_i$, $b_i^{\prime}$, $c_i$, $c_i^{\prime}$,
$1\leq i\leq 5$, and twenty-five  facets $\delta$, $\alpha_j$,
$\lambda_j$, $\gamma_j$, $1 \leq j \leq 8$ given as follows\,:

$$
\begin{array}{llll}
~~ \delta = a_1a_2b_1b_2c_2c_1, &&& \\
\alpha_1 = a_1a_2a_4b_1b_2c_2, & \alpha_2 = a_1a_2a_3a_4b_1b_2, &
\alpha_3 = a_1a_2a_3a_4a_5b_1, & \alpha_4 =
a_2a_3a_4a_5b_1c_5^{\prime}, \\
\alpha_5= a_3a_4a_5b_1 c_5^{\prime}c_4^{\prime}, & \alpha_6
=a_3a_4a_5 c_3^{\prime}c_4^{\prime}c_5^{\prime}, & \alpha_7 =
a_3a_5 c_2^{\prime} c_3^{\prime}c_4^{\prime}c_5^{\prime}, &
\alpha_8 = c_1^{\prime}c_2^{\prime}c_3^{\prime}c_4^{\prime}
c_5^{\prime} a_3, \\
\lambda_1 = a_1a_2b_2c_1c_2c_4, & \lambda_2 = a_1a_2c_1c_2c_3c_4,
& \lambda_3 = a_1c_1c_2c_3c_5c_4, & \lambda_4= a_1c_2c_3c_4c_5
b_5^{\prime}, \\
\lambda_5 = a_1c_3c_4c_5b_4^{\prime}b_5^{\prime}, & \lambda_6 =
c_3c_4c_5b_3^{\prime}b_4^{\prime}b_5^{\prime}, & \lambda_7 = c_3
c_5 b_2^{\prime}b_3^{\prime}b_4^{\prime} b_5^{\prime}, &
\lambda_8 = b_1^{\prime}b_2^{\prime}b_3^{\prime} b_4^{\prime}
b_5^{\prime}
c_3, \\
\gamma_1=a_2b_1b_2b_4c_2c_1, & \gamma_2 = b_1b_2b_3b_4c_1c_2, &
\gamma_3 = b_1b_2b_3b_4b_5c_1, & \gamma_4 = a_5^{\prime}
b_2b_3b_5b_4c_1,  \\
\gamma_5 = a_4^{\prime}a_5^{\prime}b_3b_4b_5c_1, & \gamma_6 =
a_3^{\prime}a_4^{\prime}a_5^{\prime}b_3b_5b_4, & \gamma_7 =
a_2^{\prime}a_3^{\prime}a_4^{\prime}a_5^{\prime}b_3b_5, &\gamma_8
= a_1^{\prime}a_2^{\prime}a_3^{\prime}a_4^{\prime} a_5^{\prime}
b_3.
\end{array}
$$
The dual graph $\Lambda(B^5_{30})$ is the following tree.

%%%%%%%%%%%%%%%%%%%%%%New Picture %%%%%%%%%%%

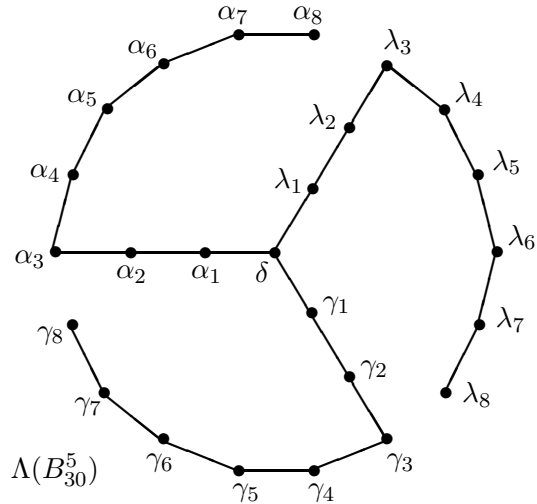
\setlength{\unitlength}{5mm}

\begin{picture}(14,14)(-10,0)

%%%%%%%%%%%%%%%%%%%%%%% %%%%%%%%%%%%%%%%%%%%%%%%%%%%%

\thicklines

\put(1.1,7){\line(1,0){5.9}} \put(7,7){\line(3,5){3}}
\put(7,7){\line(3,-5){3}}

\put(1.1,7){\line(1,4){0.5}} \put(1.6,9){\line(1,2){0.9}}
\put(2.5,10.8){\line(5,4){1.5}} \put(4,12){\line(5,2){2}}
\put(6,12.8){\line(1,0){2}}

\put(1.6,5){\line(1,-2){0.9}} \put(2.5,3.2){\line(5,-4){1.5}}
\put(4,2){\line(5,-2){2}} \put(6,1.2){\line(1,0){2}}
\put(8,1.2){\line(5,2){2}}

\put(12.9,7){\line(-1,4){0.5}} \put(12.4,9){\line(-1,2){0.9}}
\put(11.5,10.8){\line(-5,4){1.5}} \put(12.9,7){\line(-1,-4){0.5}}
\put(12.4,5){\line(-1,-2){0.9}}

\put(1,6.82){\mbox{$\bullet$}} \put(3.01,6.8){\mbox{$\bullet$}}
\put(4.98,6.8){\mbox{$\bullet$}} \put(6.84,6.8){\mbox{$\bullet$}}

\put(7.84,8.5){\mbox{$\bullet$}}
\put(8.82,10.13){\mbox{$\bullet$}}
\put(9.82,11.77){\mbox{$\bullet$}}
\put(7.82,5.2){\mbox{$\bullet$}} \put(8.82,3.5){\mbox{$\bullet$}}
\put(9.82,1.83){\mbox{$\bullet$}}

\put(7.88,1){\mbox{$\bullet$}} \put(5.88,1){\mbox{$\bullet$}}
\put(3.87,1.83){\mbox{$\bullet$}} \put(2.3,3.07){\mbox{$\bullet$}}
\put(1.45,4.88){\mbox{$\bullet$}}

\put(1.48,8.88){\mbox{$\bullet$}}
\put(3.88,11.82){\mbox{$\bullet$}}
\put(2.38,10.62){\mbox{$\bullet$}}
\put(5.88,12.6){\mbox{$\bullet$}}
\put(7.88,12.6){\mbox{$\bullet$}}

\put(12.75,6.82){\mbox{$\bullet$}}
\put(12.24,8.88){\mbox{$\bullet$}}
\put(11.35,10.59){\mbox{$\bullet$}}
\put(11.38,3.07){\mbox{$\bullet$}}
\put(12.28,4.88){\mbox{$\bullet$}}

\put(2.8,6.3){\mbox{$\alpha_2$}} \put(5.5,13.2){\mbox{$\alpha_7$}}
\put(0.5,9){\mbox{$\alpha_4$}} \put(1.5,11){\mbox{$\alpha_5$}}
\put(3.1,12.3){\mbox{$\alpha_6$}} \put(0.1,6.7){\mbox{$\alpha_3$}}
\put(7.5,13.2){\mbox{$\alpha_8$}} \put(4.8,6.3){\mbox{$\alpha_1$}}
\put(6.5,6.2){\mbox{$\delta$}}

\put(7,8.7){\mbox{$\lambda_1$}} \put(7.9,10.4){\mbox{$\lambda_2$}}
\put(9.9,12.4){\mbox{$\lambda_3$}}
\put(11.8,11){\mbox{$\lambda_4$}}
\put(12.8,9.2){\mbox{$\lambda_5$}} \put(12,3){\mbox{$\lambda_8$}}
\put(13.2,7.2){\mbox{$\lambda_6$}}
\put(12.9,5){\mbox{$\lambda_7$}}

\put(8.3,5.5){\mbox{$\gamma_1$}} \put(9.3,3.8){\mbox{$\gamma_2$}}
\put(10,1.3){\mbox{$\gamma_3$}} \put(7.9,0.5){\mbox{$\gamma_4$}}
\put(5.9,0.5){\mbox{$\gamma_5$}} \put(3.6,1.3){\mbox{$\gamma_6$}}
\put(1.7,2.8){\mbox{$\gamma_7$}} \put(0.7,4.7){\mbox{$\gamma_8$}}

\put(0,1){\mbox{$\Lambda(B^5_{30})$}}

\end{picture}

Thus, $B^5_{30}$ is a $30$-vertex stacked 5-ball, and its boundary
$S^{\,4}_{30}$ is a 30-vertex stacked 4-sphere. Let $M^4_{15}$ be
the simplicial complex obtained from $S^{\,4}_{30} \setminus
\{a_1a_2a_3a_4a_5, b_1b_2b_3b_4b_5$, $c_1c_2c_3c_4c_5,
a_1^{\prime}a_2^{\prime}a_3^{\prime}a_4^{\prime}a_5^{\prime},
b_1^{\prime}b_2^{\prime}b_3^{\prime}b_4^{\prime}b_5^{\prime},
c_1^{\prime}c_2^{\prime}c_3^{\prime}c_4^{\prime}c_5^{\prime}\}$ by
the identifications $a_i^{\prime}\equiv a_i$, $b_i^{\prime} \equiv
b_i$, $c_i^{\prime}\equiv c_i$, $1\leq i\leq 5$. It is easy to see
that each of these three identifications is admissible for
$S^{\,4}_{30}$, and remains so when the other two identifications
are already made. (Just verify that there is exactly one edge
among the four vertices $a_i$, $a_i^{\prime}$, $b_j$,
$b_j^{\prime}$ for each $i, j$, and similarly for $b$'s and $c$'s
or $c$'s and $a$'s.) Therefore, $M^4_{15}$ is indeed a
2-neighborly 4-manifold in the class ${\cal K}(4)$, with $\beta_1
= 3$ and hence $\chi = - 4$. (If $N^5_{15}$ is the simplicial
complex obtained from $B^5_{30}$ by the above identification, then
$M^4_{15}$ is the boundary of $N^5_{15}$. Then $\Lambda(N^5_{15})$
can be obtained from $\Lambda(B^5_{30})$ by adding three more
edges $\alpha_8\lambda_3$, $\lambda_8\gamma_3$ and
$\gamma_8\alpha_3$.)

Notice that the permutation $\prod_{i = 1}^5 (a_i, b_i,
c_i)(a_i^{\prime}, b_i^{\prime}, c_i^{\prime})$ is an
automorphism of order 3 in $B^5_{30}$ which induces the
automorphism $\prod_{i = 1}^5 (a_i, b_i, c_i)$ of $M^4_{15}$.
Observe that the degrees of the edges in $M^4_{15}$ are given by
the following two tables (and the above automorphism):

\begin{center}
\begin{tabular}{ccc}
\begin{tabular}{|l|ccccc|}
\multicolumn{6}{c}{Edges  of type $a_ia_j$} \\
\multicolumn{6}{c}{} \\[-1.5mm] \hline
&&&&& \\[-3.5mm]
 & $a_1$ & $a_2$ & $a_3$ & $a_4$ & $a_5$  \\[1mm]
\hline
&&&&& \\[-3.5mm]
$a_1$ & -- & 10 & 6 & 7 & 5 \\[1mm]
\hline
&&&&& \\[-3.5mm]
$a_2$ & 10 & -- & 8 & 9 & 7 \\[1mm]
\hline
&&&&& \\[-3.5mm]
$a_3$ & 6 & 8 & -- & 11 & 11 \\[1mm]
\hline
&&&&& \\[-3.5mm]
$a_4$ & 7 & 9 & 11 & -- & 11 \\[1mm]
\hline
&&&&& \\[-3.5mm]
$a_5$ & 5 & 7 & 11 & 11 & -- \\[1mm]
\hline
\end{tabular}
& &
\begin{tabular}{|l|ccccc|}
\multicolumn{6}{c}{Edges of type $a_ib_j$} \\
\multicolumn{6}{c}{} \\[-1.5mm]
\hline
&&&&& \\[-3.5mm]
 & $a_1$ & $a_2$ & $a_3$ & $a_4$ & $a_5$  \\[1mm]
\hline
&&&&& \\[-3.5mm]
$b_1$ & 7 & 9 & 7 & 8 & 6 \\[1mm]
\hline
&&&&& \\[-3.5mm]
$b_2$ & 7 & 8 & 4 & 5 & 4 \\[1mm]
\hline
&&&&& \\[-3.5mm]
$b_3$ & 4 & 5 & 6 & 7 & 8 \\[1mm]
\hline
&&&&& \\[-3.5mm]
$b_4$ & 4 & 4 & 4 & 5 & 6 \\[1mm]
\hline
&&&&& \\[-3.5mm]
$b_5$ & 5 & 4 & 5 & 6 & 7 \\[1mm]
\hline
\end{tabular}
\end{tabular}
\end{center}

These tables clearly show that the full automorphism group of
$M^4_{15}$ is of order 3. By Proposition \ref{P1}, the
face vector of $M^4_{15}$  is $(15, 105, 230, 240, 96)$. The
following is an explicit list of the 96 facets of $M^4_{15}$.

\vspace{-5mm}

\begin{eqnarray*}
&& a_1a_2b_1b_2c_1, ~ a_1b_1b_2c_1c_2, ~ a_1a_2b_1c_1c_2, ~
a_1a_2a_4b_1c_2, ~ a_2b_1b_2b_4c_1, ~ a_1b_2c_1c_2c_4, \\
&& a_1a_2a_4b_2c_2, ~ a_2b_1b_2b_4c_2, ~ a_2b_2c_1c_2c_4, ~
a_1a_4b_1b_2c_2, ~ a_2b_1b_4c_1c_2, ~ a_1a_2b_2c_1c_4, \\
&& a_1a_2b_2c_2c_4, ~ a_2a_4b_1b_2c_2, ~ a_2b_2b_4c_1c_2, ~
a_1a_2a_3a_4b_2, ~ b_1b_2b_3b_4c_2, ~ a_2c_1c_2c_3c_4, \\
&& a_1a_2a_3b_1b_2, ~ b_1b_2b_3c_1c_2, ~ a_1a_2c_1c_2c_3, ~
a_1a_2c_1c_3c_4, ~ a_1a_3a_4b_1b_2, ~ b_1b_3b_4c_1c_2, \\
&& a_1a_2c_2c_3c_4, ~ a_2a_3a_4b_1b_2, ~ b_2b_3b_4c_1c_2, ~
a_1a_2a_3a_5b_1, ~ b_1b_2b_3b_5c_1, ~ a_1c_1c_2c_3c_5, \\
&& a_1a_2a_4a_5b_1, ~ b_1b_2b_4b_5c_1, ~ a_1c_1c_2c_4c_5, ~
a_1a_3a_4a_5b_1, ~ b_1b_3b_4b_5c_1, ~ a_1c_1c_3c_4c_5, \\
&& a_2a_3a_4a_5c_5, ~ a_5b_2b_3b_4b_5, ~ b_5c_2c_3c_4c_5, ~
a_2a_3a_4b_1c_5, ~ a_5b_2b_3b_4c_1, ~ a_1b_5c_2c_3c_4, \\
&& a_2a_3a_5b_1c_5, ~ a_5b_2b_3b_5c_1, ~ a_1b_5c_2c_3c_5, ~
a_2a_4a_5b_1c_5, ~ a_5b_2b_4b_5c_1, ~ a_1b_5c_2c_4c_5, \\
&& a_3a_4a_5b_1c_4, ~ a_4b_3b_4b_5c_1, ~ a_1b_4c_3c_4c_5, ~
a_3a_4b_1c_4c_5, ~ a_4a_5b_3b_4c_1, ~ a_1b_4b_5c_3c_4, \\
&& a_3a_5b_1c_4c_5, ~ a_4a_5b_3b_5c_1, ~ a_1b_4b_5c_3c_5, ~
a_4a_5b_1c_4c_5, ~ a_4a_5b_4b_5c_1, ~ a_1b_4b_5c_4c_5, \\
&& a_3a_4a_5c_3c_4, ~ a_3a_4b_3b_4b_5, ~ b_3b_4c_3c_4c_5, ~
a_3a_4a_5c_3c_5, ~ a_3a_5b_3b_4b_5, ~ b_3b_5c_3c_4c_5, \\
&& a_3a_4c_3c_4c_5, ~ a_3a_4a_5b_3b_4, ~ b_3b_4b_5c_3c_4, ~
a_4a_5c_3c_4c_5, ~ a_3a_4a_5b_4b_5, ~ b_3b_4b_5c_4c_5, \\
&& a_3a_5c_2c_3c_4, ~ a_2a_3a_4b_3b_5, ~ b_2b_3b_4c_3c_5, ~
a_3a_5c_2c_3c_5, ~ a_2a_3a_5b_3b_5, ~ b_2b_3b_5c_3c_5, \\
&& a_3a_5c_2c_4c_5, ~ a_2a_4a_5b_3b_5, ~ b_2b_4b_5c_3c_5, ~
a_2a_3a_4a_5b_5, ~ b_2b_3b_4b_5c_5, ~ a_5c_2c_3c_4c_5, \\
&& a_1a_2a_3a_4b_3, ~ b_1b_2b_3b_4c_3, ~ a_3c_1c_2c_3c_4, ~
a_1a_2a_3a_5b_3, ~ b_1b_2b_3b_5c_3, ~ a_3c_1c_2c_3c_5, \\
&& a_1a_2a_4a_5b_3, ~ b_1b_2b_4b_5c_3, ~ a_3c_1c_2c_4c_5, ~
a_1a_3a_4a_5b_3, ~ b_1b_3b_4b_5c_3, ~ a_3c_1c_3c_4c_5.
\end{eqnarray*}

If we take the simplices $\delta$, $\alpha_1, \dots, \alpha_8$,
$\lambda_1, \dots, \lambda_8$, $\gamma_1, \dots, \gamma_8$ given
above as positively oriented simplices then that gives a coherent
orientation on $B^5_{30}$. This orientation gives a coherent
orientation on $S^4_{30}$ in which $b_5^{\prime}c_2c_3c_4c_5$,
$a_1c_3c_4c_5c_1$,  $a_1c_4c_5c_1c_2$, $a_1c_5c_1c_2c_3$,
$a_2c_1c_2c_3c_4$,
$a_5c_2^{\prime}c_3^{\prime}c_4^{\prime}c_5^{\prime}$,
$a_3c_3^{\prime}c_4^{\prime}c_5^{\prime}c_1^{\prime}$,
$a_3c_4^{\prime}c_5^{\prime}c_1^{\prime}c_2^{\prime}$,
$a_3c_5^{\prime}c_1^{\prime}c_2^{\prime}c_3^{\prime}$,
$a_3c_1^{\prime}c_2^{\prime}c_3^{\prime}c_4^{\prime}$ are
positively oriented.

Let $X = S^4_{30} \setminus \{c_1c_2c_3c_4c_5, c_1^{\prime}
c_2^{\prime}c_3^{\prime}c_4^{\prime}c_5^{\prime}\}$. Let $Y$ be
obtained from $X$ by the identifications $c_i^{\prime} \equiv
c_i$, $1\leq i\leq 5$. Then $X$ triangulates $S^3\times [0, 1]$
and the above orientation on $S^4_{30}$ induces a coherent
orientation on $X$ with positively oriented simplices (on the
boundary of $X$) $c_2c_3c_4c_5$, $c_3c_4c_5c_1$, $c_4c_5c_1c_2$,
$c_5c_1c_2c_3$, $c_1c_2c_3c_4$,
$c_2^{\prime}c_3^{\prime}c_4^{\prime}c_5^{\prime}$,
$c_3^{\prime}c_4^{\prime}c_5^{\prime}c_1^{\prime}$,
$c_4^{\prime}c_5^{\prime}c_1^{\prime}c_2^{\prime}$,
$c_5^{\prime}c_1^{\prime}c_2^{\prime}c_3^{\prime}$,
$c_1^{\prime}c_2^{\prime}c_3^{\prime}c_4^{\prime}$. This implies
that the geometric carrier of $Y$ is the twisted product $\TPSS$
(cf. \cite[pages 134--135]{st}). So, $Y$ is non-orientable. Since
$M^4_{15}$ is obtained from $Y$ by attaching two more handles, it
follows that $M^4_{15}$ is non-orientable. Therefore, Proposition
\ref{P3} implies\,:

\begin{theo}$\!\!\!${\bf .} \label{T1}
$M^4_{15}$ is a $15$-vertex triangulation of the $4$-manifold
$(\TPSS)^{\#3}$.
\end{theo}
Also, Proposition 2 implies\,:
\begin{theo}$\!\!\!${\bf .} \label{T2}
The face vector $(15, 105, 230, 240, 96)$ of $M^4_{15}$ is the
component-wise minimum over the face vectors of all triangulations
of $(\TPSS)^{\#3}$. Also, $M^4_{15}$ is a $2$-neighborly member of
Walkup's class ${\cal K}(4)$.
\end{theo}
Now, Proposition \ref{P4} implies\,:
\begin{theo}$\!\!\!${\bf .} \label{T3}
$M^4_{15}$ is a tight triangulation of $(\TPSS)^{\#3}$.
\end{theo}

It may be remarked that there are only three known infinite
families of tight combinatorial manifolds\,: (a) The standard
$d$-sphere $S^{\,d}_{d+2}$ (i.e., the boundary complex of the
simplex of dimension $d+1$), (b) the 2-neighborly 2-dimensional
triangulated manifolds (i.e., the so-called regular cases in
Heawood's map color theorem, cf. \cite{ri}) and (c) the $(2d+
3)$-vertex combinatorial $d$-manifold $K^d_{2d+3}$ due to
K\"{u}hnel (\cite{bd8, ku}): it triangulates $\TPSSD$ if $d$ is
odd and $S^{\,d-1}\!\times S^1$ if $d$ is even. Apart from these
three infinite families, only fourteen sporadic examples of tight
combinatorial manifolds were known so far, namely (i) the 9-vertex
complex projective plane due to K\"{u}hnel (\cite{kb}), (ii) three
15-vertex triangulations of homology  $\HH P^{\,2}$ due to Brehm and K\"{u}hnel
(\cite{bk}) and three more due to Lutz (\cite{lu}), (iii) a 16-vertex
triangulated K3 surface due to Casella and K\"{u}hnel (\cite{ck}),
(iv) a 15-vertex triangulation of $(\TPSS)\# (\CC P^{\,2})^{\#3}$,
(v) a 13-vertex triangulation of the homogeneous 5-manifold
$SU(3)/SO(3)$, (vi) two 12-vertex triangulations of $S^{\,3}\times
S^{\,2}$, and (vii) two 13-vertex triangulations of $S^{\,3}\times
S^{\,3}$. The last six are due to K\"{u}hnel and Lutz
(\cite{kl}). The combinatorial manifold $M^4_{15}$ constructed
here is a new entry in this select club of sporadic tight
combinatorial manifolds.

\bigskip

\noindent {\bf Acknowledgement\,:} The authors thank the anonymous
referee for many useful comments. The second author was partially
supported by DST (Grant: SR/S4/MS-272/05) and by UGC-SAP/DSA-IV.

\end{document}